\documentclass{amsart}

\usepackage{latexsym, amssymb, amsmath, amscd}
\usepackage{amsmath}

\newtheorem{proposition}{Proposition}
\newtheorem{definition}{Definition}
\newtheorem{theorem}{Theorem}
\newtheorem{lemma}{Lemma}
\newtheorem{remark}{Remark}
\newtheorem{example}{Example}
\newtheorem{corollary}{Corollary}
\newtheorem{theoremA}{Theorem}

\begin{document}

\title{Cones and Asymptotic Invariants of Multigraded Systems of Ideals}
\author{A. Wolfe}
\maketitle

\section[Introduction]{Introduction}

Recent work has discovered perhaps unexpectedly rich structures for base loci and asymptotic invariants on the cone of big divisors of smooth complex projective varieties.  One may ask what sorts of cones and functions can occur.  In general this question is not currently well-understood. However, a key feature of this work (implicit in [7] and explicit in [3]) is that the structures depend on algebraic relations between base ideals, namely, the fact that they form multigraded systems of ideals.  The question then arises: what can occur in this abstract setting?  In this paper, we give examples showing that the only restrictions obeyed by the cones and asymptotic invariants of general multigraded systems of ideals are ones of convexity, which are imposed formally.

Turning to a detailed introduction, let $X$ be a smooth, irreducible complex projective variety, and $N^1(X)_{\mathbf{R}}$ be the finite-dimensional vector space of real numerical equivalence classes of divisors on $X$.  This vector space contains interesting cones reflecting the geometry of $X$.  For example, one has the cone of nef divisors Nef$(X) \subset N^1(X)_{\mathbf{R}}$, the closure of the cone spanned by ample divisors on $X$, and $\overline{\textnormal{Eff}(X)} \subset N^1(X)_{\mathbf{R}}$ the closure of the cone spanned by effective divisors, which contains the cone of big divisors Big$(X)$ as its interior.  These cones have attracted a great deal of attention and provide many examples of interesting behavior.

Moreover, it has recently been realized that there are naturally defined functions on these cones that reflect the behavior of the linear series in question.  For example, fixing $x \in X$ and an effective divisor $D$ on $X$, set
$$
\mbox{ord}_x ( |D| ) := \min_{D' \in |D|} \mbox{mult}_x (D')
$$
$$
\mbox{ord}_x ( ||D||) := \lim_{m \rightarrow \infty} \frac{ \mbox{ord}_x (|mD|)}{m}
$$
Thus $\textnormal{ord}_x (|D|)$ measures the multiplicity of a general divisor in the linear series $|D|$ and $\textnormal{ord}_x (|| D ||)$ reflects its asymptotic growth.  Nakayama proves that this depends only on the numerical equivalence class of the divisor $D$ and extends by continuity to a function 
$$
\mbox{ord}_x : \textnormal{Big}(X) \rightarrow \mathbf{R}
$$

Precisely what sort of cones and functions occur remains a tantalizing open problem.  However, a key observation is that the cones and functions do not depend on the global geometry of the variety $X$, but but can be defined in a purely algebraic setting, in terms of ideals.  If one chooses divisors $D_1, \ldots , D_{\rho}$ whose classes span $N^1 (X)_{\mathbf{R}}$, the nef cone will be given by the closure of the cone spanned by the points $(a_1, \ldots a_{\rho})$ for which the corresponding the base ideals
$$
\mathbf{b} (| ma_1 D_1 + \cdots + ma_{\rho} D_{\rho} |) = \mathcal{O}_X \textnormal { for } m \gg 0.
$$
The pseudoeffective cone has a similar interpretation.

We therefore are led to consider formal systems of ideals, taking as an axiom the rule
$$
\mathbf{b}(|{D_1}|) \mathbf{b}(|{D_2}|) \subseteq \mathbf{b}(|{D_1+D_2}|)
$$
arising from the map
$$
H^0 (X, \mathcal{O}(D_1) ) \otimes H^0 (X, \mathcal{O}(D_2)) \rightarrow H^0 (X, \mathcal{O}(D_1 + D_2)).
$$
Besides providing a method of constructing examples of behavior that can occur globally on projective varieties, this setting captures sufficient information to define analogues of the cones and asymptotic functions that have been discovered in the projective setting.

\begin{definition}  Let $X$ be a smooth, irreducible complex variety.  A $\mathbf{Z}^{\rho}$-graded system of ideals on $X$ is a collection of ideal sheaves $\mathbf{a}_{\bullet} = \{ \mathbf{a}_v \}_{v \in \mathbf{Z}^{\rho}}$ so that for $v, w \in \mathbf{Z}^{\rho}$, one has the relation
$$
\mathbf{a}_v \mathbf{a}_w \subseteq \mathbf{a}_{v+w}.
$$
\end{definition}

We note that the essential examples in this paper occur on affine varieties, so the reader may safely take $X = \textnormal{Spec}(R)$ to be affine and $\mathbf{a}_v$ to be ideals of $R$.  In fact, there is no loss in taking $R$ to be a polynomial ring.

One defines the algebraic Neron-Severi space $\textnormal{NS}_{\mathbf{R}}(\mathbf{a}_{\bullet})$ of a graded system $\mathbf{a}_{\bullet}$  to be $\mathbf{Z}^{\rho} \otimes \mathbf{R}$, the analog of $N_{\mathbf{R}}^1(X)$.  Writing $\langle S \rangle$ for the closed convex cone spanned by a set $S \subseteq \mathbf{Z}^{\rho} \otimes \mathbf{R}$, one has within $\textnormal{NS}_{\mathbf{R}}(\mathbf{a}_{\bullet})$ the cones
$$
\textnormal{Nef} (\mathbf{a}_{\bullet}) = \langle \{ v \in \mathbf{Z}^{\rho} | \mathbf{a}_v = (1) \} \rangle .
$$
$$
\overline{\textnormal{Eff}}(\mathbf{a}_{\bullet}) = \langle \{ v \in \mathbf{Z}^{\rho} | \mathbf{a}_v \neq 0 \} \rangle .
$$
Again by analogy to the projective case, one may call the interior of these cones the ample and big cones of the graded system, respectively.

All the naturally defined asymptotic invariants in the projective setting have their analogs in the algebraic setting.  For example, given a prime ideal $\mathbf{p} \subset R$ one sets
$$
\textnormal{ord}_{\mathbf{p}} (\mathbf{a}_{\bullet}, v) := \limsup_n \frac{ord_{\mathbf{p}}\mathbf{a}_{nv}}{n}
$$
In section 3 we will define two other asymptotic invariants, the log-canonical threshold and the multiplicity (for systems of zero-codimensional ideals).

It is natural to ask how the asymptotic functions associated to these cones behave in this general setting.  It is hoped that many of these invariants are very well-behaved on projective varieties.  For example, in [1] Campana and Peternell prove a real version of the Nakai-Moishezon criterion to show that the nef cone of a projective variety is locally given at almost all points by a polynomial inequality.  One could also hope that the asymptotic invariants described above are also locally analytic.  However, our examples show that in the abstract setting of systems of ideals, the resulting invariants can be fairly wild.

\begin{theoremA}
Given a closed convex cone $C$ contained in a finite-dimensional real vector space $V$, with nonempty interior, and with the vertex of $C$ lying at the origin of $V$, there is a graded system of ideals having $C$ as its nef cone.
\end{theoremA}
It follows quickly that given a pair of closed convex cones $C_1 \subseteq C_2$, there exists a graded system of ideals having nef cone $C_1$ and effective cone $C_2$.  Moreover, the construction used to prove the theorem shows that there exist graded systems whose asymptotic invariants are not differentiable on any open set in $\textnormal{Big} (\mathbf{a}_{\bullet})$.

The following theorem shows moreover that the regularity of the asymptotic invariants does not follow formally from the regularity of the nef or effective cones of a graded system of ideals.  
\begin{theoremA}
There exists a graded system $\{ \mathbf{a}_v \}_{v \in \mathbf{Z}^2}$ so that the boundary of $\textnormal{Nef}(\mathbf{a}_{\bullet})$ is given by linear functions, but for which $\textnormal{ord}_0$ is nowhere differentiable on an open set in $\textnormal{NS}_{\mathbf{R}}(\mathbf{a}_{\bullet})$.
\end{theoremA}

\subsection{Acknowledgments}
This paper was prepared under the direction of Robert Lazarsfeld.  We are grateful for his suggestions and encouragement.  We are indebted to Jason Howald for the idea of using graded systems of monomial ideals defined by convex regions with complicated boundaries.  We thank Alexander Barvinok, Juha Heinonen, Mattias Jonsson, and Dror Varolin for their help and useful comments.  We also thank the University of Michigan Mathematics Department and the National Science Foundation's VIGRE program for their support while this article was in preparation.

\section[Multigraded Systems of Ideal Sheaves]{Multigraded Systems of Ideal Sheaves}

We give the definition of multigraded systems of ideal sheaves and illustrate the formalism with a number of examples.

\begin{definition}
Let $X$ be a smooth, irreducible complex variety.  A $\mathbf{Z}^{\rho}$-graded system of ideal sheaves on $X$ is a collection of ideal sheaves
$ \mathbf{a}_{\bullet} = \{ \mathbf{a}_v \}_{v \in \mathbf{Z}^{\rho}}$ so that, for $v,w \in \mathbf{Z}^{\rho}$, one has the relation
$$
\mathbf{a}_v \mathbf{a}_w \subseteq \mathbf{a}_{v+w}
$$
\end{definition}

\begin{remark}  \textnormal{We will often give examples of graded systems of ideal sheaves where the indexing set is some subsemigroup } $S$ \textnormal{ of } $\mathbf{Z}^{\rho}$, \textnormal{ most often } $\mathbf{N}^{\rho} \subset \mathbf{Z}^{\rho}$.  \textnormal{ In fact, this is equivalent to giving a } $\mathbf{Z}^{\rho}$-\textnormal{graded system, since one may recover a } $\mathbf{Z}^{\rho}$-\textnormal{graded system by setting the ideal sheaves } $\mathbf{a}_{v} := 0$ \textnormal{ for } $v \in \mathbf{Z}^{\rho} \setminus S$.
\end{remark}

Conversely, given a $\mathbf{Z}^{\rho}$-graded system $\mathbf{a}_{\bullet}$ and a subsemigroup $S \subseteq \mathbf{Z}^{\rho}$, one may consider the $\mathbf{Z}^{\rho}$-graded system $\mathbf{b}_{\bullet}$, the truncation of $\mathbf{a}_{\bullet}$ by $S$, defined by
$$
\mathbf{b}_v = \mathbf{a}_v, v \in S
$$
$$
\mathbf{b}_v = 0, v \notin S.
$$

We will soon exhibit multigraded systems of ideal sheaves with effective cones equal to the entire algebraic Neron-Severi space.  We will then be able to exhibit a graded system having effective cone $C$, for an arbitrary closed convex cone $C$, by taking the truncation of a graded system with $\overline{\textnormal{Eff}}(\mathbf{a}_{\bullet}) = \textnormal{NS}(\mathbf{a}_{\bullet})$ by the subsemigroup $C \cap \mathbf{Z}^{\rho}$ of $\mathbf{Z}^{\rho}$.

\begin{remark}  \textnormal{We will usually consider only cases where } $\overline{\textnormal{Eff}}(X)$ \textnormal{has nonempty interior inside } $\mathbf{Z}^{\rho} \otimes \mathbf{R}$.
\end{remark}

\begin{example} \textnormal{The motivating geometric example is the following: let} $D_1, \ldots, D_{\rho}$ \textnormal{be divisors on a smooth projective variety} $X$.  \textnormal{Set}
$$
\mathbf{a}_{(n_1, \ldots, n_k)} := \mathbf{b} (|n_1 D_1 + \cdots + n_{\rho} D_{\rho}|).
$$
\end{example}

\begin{example} \textnormal{  This example should be considered the trivial case: choose ideal sheaves} $I_1, \ldots , I_{\rho}$ \textnormal{and set}
$$
\mathbf{a}_{(n_1, \ldots, n_{\rho})} := I_1^{n_1} \cdots I_{\rho}^{n_{\rho}}.
$$
\textnormal{Here} $I^n = (1)$ \textnormal{if} $n \leq 0$.
\end{example}

\begin{example}
\textnormal{Colon ideals give another example of naturally occurring graded systems.  Fix} $X = \textnormal{Spec}(R)$ \textnormal{a smooth irreducible affine variety.  Given a graded system of ideal}s $ \{ \mathbf{a}_v \}_{v \in \mathbf{Z}^{\rho}}$ \textnormal{and an ideal} $I \subset R$, \textnormal{one may form a} $\mathbf{Z}^{\rho} \times \mathbf{N}$-\textnormal{graded system of ideals by setting}
$$
\mathbf{b}_{(m,n)}:= (\mathbf{a}_m : I^n)
$$
\end{example}

\begin{example}
\textnormal{Fix a proper, birational map} $\mu : X' \rightarrow X$ \textnormal{with exceptional divisors} $D_1, \ldots , D_{\rho}$.  \textnormal{For example, one may take} $X = \textnormal{Spec}(R)$ \textnormal{ and } $X'$ \textnormal{a resolution of singularities of the blow-up of an ideal} $I \subset R$.  \textnormal{Thanks to the natural map}
$$
\mu_{*} \mathcal{O}_{X'}(\sum_{i=1}^{\rho} -n_i D_i) \otimes \mu_{*} \mathcal{O}_{X'}(\sum_{i=1}^{\rho} -m_i D_i) \rightarrow \mu_{*} \mathcal{O}_{X'}(\sum_{i=1}^{\rho} -(n_i+m_i) D_i)
$$
\textnormal{one obtains an} $\mathbf{N}^{\rho}$-\textnormal{graded system of ideal sheaves on} $X$ \textnormal{by setting} 
$$
\mathbf{a}_{(n_1, \ldots n_{\rho})} = \mu_{*} \mathcal{O}_{X'}(\sum_{i=1}^{\rho} -n_i D_i).
$$
\textnormal{As a special case, take } $X= \textnormal{Spec}(R)$, $I$ \textnormal{a prime ideal, and} $X'$ \textnormal{a smooth variety dominating} $\textnormal{Bl}_I (X)$, \textnormal{the blow-up of} $X$ \textnormal{at} $I$.  \textnormal{Let} $E$ \textnormal{be the unique prime exceptional divisor in} $\textnormal{Bl}_I (X)$ \textnormal{dominating the smooth locus of} $V(I)$.  \textnormal{Then if one writes} $E'$ \textnormal{for the proper transform of} $E$ \textnormal{on} $X'$, \textnormal{one has that the symbolic powers} $I^{ \langle n \rangle}$ \textnormal{are given by}
$$
I^{ \langle n \rangle} = \mu_{*} \mathcal{O}_{X'} (-nE').
$$
\end{example}

\begin{example}
\textnormal{ Fix a closed, convex region} $P \subset \mathbf{R}_{\geq 0}^n$ \textnormal{with the property that if} $p \in P$, \textnormal{then} $p + \mathbf{R}_{\geq 0}^n \subseteq P$.  \textnormal{This guarantees that} $ \{ x^v | v \in P \} $ \textnormal{forms a basis over} $\mathbf{C}$ \textnormal{for a monomial ideal} $\mathbf{a} \subset R = \mathbf{C} [ x_1, \ldots x_n]$.  \textnormal{We now define an} $\mathbf{N}$-\textnormal{graded system of ideals on} $\mathbf{A}_{\mathbf{C}}^n$ \textnormal{by}
$$
\mathbf{a}_n := \textnormal{span}_{\mathbf{C}} \{ x^v | v \in nP \}.
$$
\textnormal{Since} $\mathbf{a}_1 \supset \mathbf{a}_2 \supset \cdots$, \textnormal{one may extend this to a} $\mathbf{Z}$-\textnormal{graded system of ideals by setting} $\mathbf{a}_n := R $ \textnormal{for all} $n \leq 0$.
\end{example}

\begin{example}  
\textnormal{We discuss a few basic constructions with graded systems of ideals.  If} $\mathbf{a}_{\bullet}$ \textnormal{and} $\mathbf{b}_{\bullet}$ \textnormal{are two} $\mathbf{Z}^{\rho}$-\textnormal{graded systems of ideals then we may construct their product and intersection:}
$$
(\mathbf{a}_{\bullet} \mathbf{b}_{\bullet})_{v} := \mathbf{a}_v \mathbf{b}_v
$$
$$
(\mathbf{a}_{\bullet} \cap \mathbf{b}_{\bullet})_v := \mathbf{a}_v \cap \mathbf{b}_v
$$
\textnormal{Also, given a homomorphism} $\phi : \mathbf{Z}^n \rightarrow \mathbf{Z}^m$ \textnormal{and a graded system} $\{ \mathbf{a}_v \} _{v \in \mathbf{Z}^m}$ \textnormal{one may define a} $\mathbf{Z}^n$-\textnormal{graded system} $\mathbf{a}_{\bullet} = \phi^{*}(\mathbf{b}_{\bullet}) $ \textnormal{by, for} $w \in \mathbf{Z}^n$, 
$$
\mathbf{a}_w := \mathbf{b}_{\phi(w)}
$$
\end{example}

\begin{example}
\textnormal{Example 6 allows us to introduce multigraded systems of monomial ideals arising from the graded systems of monomial ideals discussed in Example 5.  Let} $P_1, P_2$ \textnormal{be regions in} $\mathbf{R}_{\geq 0}^n$ \textnormal{as in Example 5, and let} $\mathbf{a}_{\bullet}^1$ \textnormal{and} $\mathbf{a}_{\bullet}^2$ \textnormal{be the graded systems defined by} 
$$
\mathbf{a}_{n}^i := \textnormal{sp}_{\mathbf{C}} \{ x^v | v \in nP_i \}.
$$
\textnormal{From these we construct two $\mathbf{Z}^2$-graded systems of ideals.  Write} $\textnormal{pr}_1, \textnormal{pr}_2 : \mathbf{Z}^2 \rightarrow \mathbf{Z}$ \textnormal{ for the first and second projections.  We define }
$$
\mathbf{b}_{\bullet} := \textnormal{pr}_1^{*} \mathbf{a}_{\bullet}^1 \cap \textnormal{pr}_2^{*} \mathbf{a}_{\bullet}^2
$$
$$
\mathbf{c}_{\bullet} := \textnormal{pr}_1^{*} \mathbf{a}_{\bullet}^1 \cdot \textnormal{pr}_2^{*} \mathbf{a}_{\bullet}^2
$$
\textnormal{For} $n,m \in \mathbf{N}$ \textnormal{the graded systems above are given by:}
$$
\mathbf{b}_{(n,m)} = \{ x^v | v \in nP \cap mQ \}
$$
$$
\mathbf{c}_{(n,m)} = \{ x^v | v \in nP + mQ \}
$$
\end{example}
We will use in particular the geometry of $P$ and $Q$ and the construction of the graded system $\mathbf{b}_{\bullet}$ to give examples of pathological behavior of asymptotic invariants of graded systems of ideals.

\section[Asymptotic Invariants of Multigraded Systems]{Asymptotic Invariants of Multigraded Systems}

The asymptotic invariants defined for base loci of line bundles on projective varities also have their natural analogs in the case of general multigraded systems of ideal sheaves.  We first recall definitions of these invariants for single ideals or ideal sheaves.

Let $X$ be a smooth, irreducible complex variety and $I$ an ideal sheaf on $X$.  For $Z \subset X$ a proper subvariety, we set 
$$
\textnormal{ord}_Z (I) := \textnormal{max } \{ n | I \subseteq I_Z^{ \langle n \rangle } \}.
$$
We recall that this is the maximum $n$ so that if $D$ is a differential operator of order at most $n-1$, $Df(p) = 0$ at all $p \in Z$ (or equivalently, at a general point $p \in Z$).

Let $x \in X$ be a point and $I$ an ideal sheaf on $X$ with $x \in \textnormal{Supp}(\mathcal{O}_X / I)$ so that $I \mathcal{O}_{X,x}$ is zero-codimensional and write $d = \textnormal{dim}(\mathcal{O}_{X,x})$.  We recall the Samuel multiplicity, defined as
$$
e_x (I) := \lim_n \frac{ l (\mathcal{O}_{X,x}/I^n)}{n^d/d!}.
$$

Finally, if $\lambda \in \mathbf{R}_{> 0}$ we may form the multiplier ideal
$\mathcal{J} (X, \lambda \cdot I)$.  We define the log-canonical threshold of $I$ to be
$$
c(I) := \textnormal{inf } \{ \lambda \in \mathbf{R}_{ > 0} | \mathcal{J} (X, \lambda \cdot I) \neq \mathcal{O}_{X} \}.
$$
We refer to [5] for definitions and basic facts about multiplier ideals.  Formally, however, the following lemma in [3] shows that it is preferable to work with the reciprocal Arn$(I) = \frac{1}{c(I)}$.

\begin{lemma}  Let $\mathbf{a}_{\bullet}$ be an $\mathbf{N}$-graded system of ideals on a smooth, irreducible complex variety $X$ of dimension $d$ and fix indices $p,q$ so that $\mathbf{a}_p, \mathbf{a}_q \neq 0$.  Then for $x \in X$ a closed point and $Z \subset X$ a proper subvariety,
$$
\textnormal{ord}_Z (\mathbf{a}_{p + q}) \leq \textnormal{ord}_Z (\mathbf{a}_p)+ \textnormal{ord}_Z (\mathbf{a}_q)
$$
$$
\textnormal{Arn} (\mathbf{a}_{p + q}) \leq \textnormal{Arn}(\mathbf{a}_p)+ \textnormal{Arn} (\mathbf{a}_q)
$$
$$
e_x(\mathbf{a}_{p+q})^{1/d} \leq e_x(\mathbf{a}_{p})^{1/d} + e_x(\mathbf{a}_{p})^{1/d}
$$
\end{lemma}

By Lemma 1.4 in [5], this implies that one can define:
$$
\textnormal{ord}_Z (\mathbf{a}_{\bullet}) := \lim_{n \rightarrow \infty} \frac{ \textnormal{ord}_Z (\mathbf{a}_n)}{n}
$$
$$
\textnormal{Arn}(\mathbf{a}_{\bullet}) :=  \lim_{n \rightarrow \infty} \frac{ \textnormal{Arn}(\mathbf{a}_n)}{n}
$$
$$
e_x (\mathbf{a}_\bullet) := \lim_{n \rightarrow \infty} \frac{ \textnormal e_x(\mathbf{a}_n)}{n^d/d!}.
$$

We apply this to the multigraded case of a $\mathbf{Z}^{\rho}$-graded system of ideals $\mathbf{a}_{\bullet}$ by fixing a vector $v \in \mathbf{Z}^{\rho}$ and examining the $\mathbf{N}$-graded system $\mathbf{b}_{\bullet}$ given by
$$
\mathbf{b}_n := \mathbf{a}_{nv}.
$$
\begin{definition}  
Let $\mathbf{a}_{\bullet}$ be a $\mathbf{Z}^{\rho}$-graded system of ideals on a smooth, irreducible complex variety $X$, and fix $Z \subset X$ a proper subvariety of $X$.  Then for $v \in \mathbf{Z}^{\rho}$, 
$$
\textnormal{ord}_Z (\mathbf{a}_{\bullet}, v) = \textnormal{ord}_Z (\mathbf{b}_{\bullet}) = \lim_{n \rightarrow \infty} \frac{ \textnormal{ord}_Z (\mathbf{a}_{nv})}{n}
$$
$$
\textnormal{Arn} (\mathbf{a}_{\bullet}, v) = \textnormal{Arn} (\mathbf{b}_{\bullet}) = \lim_{n \rightarrow \infty} \frac{ \textnormal{Arn} (\mathbf{a}_{nv})}{n}.
$$
If, moreover, $\mathbf{a}_{nv}$ is supported at $x \in X$ for all $n \gg 0$, we set 
$$
e_x (\mathbf{a}_{\bullet}, v) = e_x (\mathbf{b}_{\bullet}) = \lim_{n \rightarrow \infty} \frac{ e_x(\mathbf{a}_{nv})}{n^d/d!}
$$
\end{definition}

When it is clear which multigraded system is intended, we will suppress its specification and write only $\textnormal{ord}_Z(v)$, Arn$(v)$, or $e_x (v)$.

Since these functions are all homogeneous with respect to multiplication of $v$ by a positive integer multiple, they are well-defined on rational classes in $\mathbf{Z}^{\rho} \otimes \mathbf{Q}$.  Moreover, they are all convex as functions of $v$, since the $n^{\textnormal{th}}$ term in the sequences defining them is also convex as a function of $v$.

Recall that $\textnormal{Big}(\mathbf{a}_{\bullet})$ is defined to be the interior of $\overline{\textnormal{Eff}} (\mathbf{a}_{\bullet}) \subseteq \textnormal{NS}_{\mathbf{R}} (\mathbf{a}_{\bullet} )$.  In [3] it is shown that, when the ample cone of the multigraded system contains a basis for $\textnormal{NS}_{\mathbf{R}}(\mathbf{a}_{\bullet})$, these functions satisfy a Holder-type inequality which guarantees that the invariants extend to continuous, convex functions on the real-valued points of $\textnormal{Big} (\mathbf{a}_{\bullet})$.

\begin{remark}
\textnormal{ It is a theorem of Mustata that the invariant } $e(\mathbf{a}_{\bullet})$ \textnormal{ agrees with the volume of the system of zero-dimensional ideals} $\mathbf{a}_{\bullet}$ \textnormal{given by }
$$
\textnormal{vol}(\mathbf{a}_{\bullet}) = \limsup_{n} \frac{ (\textnormal{dim}(R))! \cdot l(R/\mathbf{a}_{nv})}{n^{\textnormal{dim}(R)}}. 
$$
\end{remark}

\section[Graded Systems of Monomial Ideals]{Graded Systems of Monomial Ideals}
The ease of calculation of invariants afforded by systems of monomial ideals on $X = \mathbf{C}^k$ makes their study attractive and gives a good source of examples.  We will focus in this section on $\mathbf{N}$-graded systems of monomial ideals.  Since, as was indicated in the discussion preceding Definition 3, asymptotic invariants for multigraded systems are calculated by first passing to $\mathbf{N}$-graded systems, there is no loss of generality when calculating asymptotic invariants.

We recall the Newton polyhedron of a monomial ideal, from which our convex sets will arise: given a monomial ideal $\mathbf{a}$, one defines its Newton polyhedron $P(\mathbf{a}) \subset \mathbf{R}_{\geq 0}^k$ to be the convex hull of the set of all exponent vectors of monomials in $\mathbf{a}$.  Since $\mathbf{a}$ is an ideal, the Newton polyhedron $P( \mathbf{a} )$ absorbs $\mathbf{R}_{\geq 0}^k$ under addition: if $v \in P$, $v + \mathbf{R}_{\geq 0}^k \subseteq P$.

We will show that the asymptotic invariants of an $\mathbf{N}-$graded system of monomial ideals $\mathbf{a}_{\bullet}$ are determined by the geometry of the convex set  
$$
P(\mathbf{a}_{\bullet}) := \cup_{n \in \mathbf{N} } \tfrac{1}{n} P(\mathbf{a}_{n}).
$$
Such sets were considered by Mustata (see [5]).  

An important point is that the sets $P(\mathbf{a}_{\bullet})$ may be fairly wild.  We observe in particular that even though each $P(\mathbf{a}_n)$ is polyhedral, $P(\mathbf{a}_{\bullet})$ need not be.  In fact, after computing the asymptotic invariants of $\mathbf{a}_{\bullet}$ in terms of $P(\mathbf{a}_{\bullet})$, we will show that, up to closure, any closed convex set in $\mathbf{R}_{ \geq 0}^k$ absorbing  $\mathbf{R}_{ \geq 0}^k$ under addition can be realized as $P(\mathbf{a}_{\bullet})$ for some graded system $\mathbf{a}_{\bullet}$.  In particular, the epigraph $P$ of any everywhere decreasing convex function on $\mathbf{R}_{>0}$ can occur as such a set.

\begin{lemma}
Let $\mathbf{a}_{\bullet}$ be an $\mathbf{N}$-graded system of ideals.  Then $P(\mathbf{a}_{\bullet}) \subseteq \mathbf{R}_{\geq 0}^k$ is a convex set absorbing $\mathbf{R}_{\geq 0}^k$ under addition.
\end{lemma}
\textbf{Proof}  If $w \in P(\mathbf{a}_{\bullet})$, then $w\in \tfrac{1}{n} P(\mathbf{a}_n)$ for some $n$, and 
$$
w + \mathbf{R}_{\geq 0}^k \subseteq \tfrac{1}{n} P(\mathbf{a}_n) \subseteq P(\mathbf{a}_{\bullet}).
$$
To show convexity, suppose $x,y \in P(\mathbf{a}_{\bullet})$ and $0 \leq \lambda \leq 1$; say $x \in \frac{1}{n} P(\mathbf{a}_{n})$, $y \in \frac{1}{m} P(\mathbf{a}_{m})$.  Then $nx \in P(\mathbf{a}_{n})$ and $my \in P(\mathbf{a}_{m})$, so $mnx, nmy \in P(\mathbf{a}_{nm})$.  By convexity of Newton polyhedra,
$$
\lambda (nmx) + (1 - \lambda)(nmy) \in P(\mathbf{a}_{mnv}),
$$
$$
\lambda x + (1 - \lambda)y \in \tfrac{1}{mn}  P(\mathbf{a}_{mnv}) \subseteq P(\mathbf{a}_{\bullet}).  \Box
$$

We now compute the asymptotic invariants in terms of $P(\mathbf{a}_{\bullet})$.  We first fix notation.  For $v = ( v_1, \ldots v_k) \in \mathbf{R}_{ \geq 0}^k$, set $| v | := \sum_{i=1}^k v_i.$  We define $\mathbf{1}$ to be the vector $(1, \ldots, 1)  \in \mathbf{R}_{ \geq 0 }^k$.  Given a region $P \subset \mathbf{R}_{ \geq 0 }^k$ absorbing $ \mathbf{R}_{ \geq 0 }^k$ under addition, we set
$$
\lambda(P) := \inf \{ \lambda \in \mathbf{R}_{\geq 0} | \lambda \cdot \mathbf{1} \in P \}.
$$

\begin{proposition}  Let $\mathbf{a}_{\bullet}$ be an $\mathbf{N}$-graded system of monomial ideals on $\mathbf{C}^k$.  Then:
\begin{enumerate}
  \item $\textnormal{ord}_0 (\mathbf{a}_{\bullet}) = \inf \Big \{ | v | : v \in P( \mathbf{a}_{\bullet}) \Big \}$. \label{ord}
  \item $\textnormal{Arn}(\mathbf{a}_{\bullet}) = \lambda \Big ( P ( \mathbf{a}_{\bullet}) \Big )$.  \label{c}
    \item $e(\mathbf{a}_{\bullet}) = k! \cdot \textnormal{Vol} \Big( \mathbf{R}_{\geq 0}^k \setminus P(\mathbf{a}_{\bullet}) \Big)$, if $\mathbf{a}_{n}$ is zero-dimensional for all $n \gg 0$. \label{e}
\end{enumerate}
\end{proposition}

We first establish two simple lemmas concerning the ideals $\mathbf{a}_{n!}$ that will simplify the proof.
\begin{lemma}
Let $\{ \mathbf{a}_n \}_{n \in \mathbf{N}}$ be a graded system of ideals.  Then 
$$
\bigcup_{n=1}^L \tfrac{1}{n} P(\mathbf{a}_{n}) \subseteq \tfrac{1}{L!} P(\mathbf{a}_{L!}).
$$
\end{lemma}
\textbf{Proof of Lemma 3}  Observe that if $m,n$ are integers with $q = \frac{n}{m} \in \mathbf{N}$, then $\mathbf{a}_{mv}^q \subseteq \mathbf{a}_{nv}$, so
$$
\tfrac{1}{m}P(\mathbf{a}_{m}) = \tfrac{1}{n} qP(\mathbf{a}_{m}) = \tfrac{1}{n} P(\mathbf{a}_{m}^q) \subseteq \tfrac{1}{n} P(\mathbf{a}_{n}).
$$
Setting $n=L!$ and $1 \leq m \leq L$, the lemma follows. $\Box$

\begin{lemma}
The sequences
$$
\frac{\textnormal{ord}_0 (\mathbf{a}_{L!})}{L!}, 
\frac{\textnormal{Arn}   (\mathbf{a}_{L!})}{L!},
\frac{e(\mathbf{a}_{L!})}{(L!)^k}
$$
are monotone decreasing sequences converging to $\textnormal{ord}_0 (\mathbf{a}_{\bullet}), \textnormal{ Arn} (\mathbf{a}_{\bullet}), e(\mathbf{a}_{\bullet})$ respectively.
\end{lemma}
\textbf{Proof of Lemma 4}  It follows from Definition 3 that the sequences converge to the stated limits.  We need only check that the sequences are monotone decreasing.  For the invariant $\textnormal{ord}_0$, since $ \mathbf{a}_{L!}^{L+1} \subseteq \mathbf{a}_{(L+1)!}$,
\begin{eqnarray*}
\frac{ \textnormal{ord}_0 (\mathbf{a}_{(L+1)!}) }{(L+1)!} &\leq& 
\frac{ \textnormal{ord}_0 (\mathbf{a}_{L!}^{L+1})}{(L+1)!} \\
&=& \frac{(L+1) \textnormal{ord}_0 (\mathbf{a}_{L!})}{(L+1)!} \\
&=& \frac{\textnormal{ord}_0 (\mathbf{a}_{L!})}{L!}.
\end{eqnarray*}
The invariants $\textnormal{Arn}( \cdot )$ and $e^{1/k} (\cdot)$, which obey formal properties similar to those of $\textnormal{ord}_0$, are handled in a similar fashion.  $\Box$

\medskip
\textbf{Proof of the Proposition}  To begin with, recall that if $\mathbf{a}$ is a single monomial ideal, 
\begin{eqnarray*}
\textnormal{ord}(\mathbf{a}) &=& \inf \Big \{ | v | : v \in P(\mathbf{a}) \Big \} \\
\textnormal{ Arn}(\mathbf{a}) &=& \lambda \Big ( P(\mathbf{a}) \Big ) \\
e(\mathbf{a}) &=& k! \cdot \textnormal{Vol} \Big ( \mathbf{R}_{ \geq 0 }^k \setminus P(\mathbf{a}) \Big ).
\end{eqnarray*}
Applying these formulas and recalling that $\tfrac{1}{n} P (\mathbf{a}_n ) \subset P( \mathbf{a}_{\bullet})$, one has
\begin{eqnarray*}
\frac{ \textnormal{ord}_0 (\mathbf{a}_n) }{n} &=& \inf \Big \{ |v| : v \in \tfrac{1}{n} P(\mathbf{a}_n) \Big \} \textnormal{ } \geq \textnormal{ } \inf \Big \{ |v| : v \in P(\mathbf{a}_{\bullet}) \Big \} \\
\frac{\textnormal{Arn}(\mathbf{a}_n)}{n} &=& \lambda \Big ( \tfrac{1}{n} P (\mathbf{a}_n ) \Big) \textnormal{ } \geq \textnormal{ } \lambda \Big ( P( \mathbf{a}_{\bullet} )   \Big ) \\
\frac{ e (\mathbf{a}_n) }{n^k} &=&  k! \cdot \textnormal{Vol} \Big (  \mathbf{R}_{ \geq 0}^k \setminus \tfrac{1}{n} P(\mathbf{a}_n) \Big ) \textnormal{ } \geq \textnormal{ }  k! \cdot \textnormal{Vol} \Big (  \mathbf{R}_{ \geq 0}^k \setminus P(\mathbf{a}_{\bullet}) \Big ) .
\end{eqnarray*}
Taking limits as $n \rightarrow \infty$ in the preceding inequalities gives
\begin{eqnarray*}
\textnormal{ord}_0 (\mathbf{a}_{\bullet}) &\geq&  \inf \Big \{ | v | : v \in P( \mathbf{a}_{\bullet}) \Big \} \\
\textnormal{Arn}(\mathbf{a}_{\bullet}) &\geq& \lambda \Big ( P ( \mathbf{a}_{\bullet}) \Big ) \\
e(\mathbf{a}_{\bullet}) &\geq& k! \cdot \textnormal{Vol} \Big( \mathbf{R}_{\geq 0}^k \setminus P(\mathbf{a}_{\bullet}) \Big).
\end{eqnarray*}

We now prove the reverse inequalities.  For (1), let $M = \textnormal{inf } \big \{ |m| : m \in P(\mathbf{a}_{\bullet}) \big \}$.  Let $\epsilon > 0$ and select a rational point $m \in P(\mathbf{a}_{\bullet})$ so that $|m| < M + \epsilon$.  Then $m \in \tfrac{1}{L} P(\mathbf{a}_L)$ for some $L$.  A fortiori $m \in \tfrac{1}{L!} P(\mathbf{a}_{L!})$ and thus for all $l \geq L$,
$$
\frac{ \textnormal{ord}_0 (\mathbf{a}_{l!})}{l!} \textnormal{ } \leq \textnormal{ } \frac{ \textnormal{ord}_0 (\mathbf{a}_{L!})}{L!} \textnormal{ } \leq \textnormal{ } |m|,
$$
But then
$$
\textnormal{ord}_0 (\mathbf{a}_{\bullet}) \textnormal{ } \leq \textnormal{ } |m| \textnormal{ } < \textnormal{ } M + \epsilon.
$$

Turning to (2), given $\epsilon > 0$, choose  $v = \beta \cdot \mathbf{1} \in P(\mathbf{a}_{\bullet})$ with $\beta < \lambda + \epsilon$.  As above $p \in \tfrac{1}{L} P(\mathbf{a}_L)$ for some $L$, and for all $l \geq L$, 
$$
\frac{\textnormal{Arn} (\mathbf{a}_{l!}) }{l!} \textnormal{ } \leq  \textnormal{ } \frac{\textnormal{Arn} (\mathbf{a}_{L!}) }{L!}  \textnormal{ } \leq  \textnormal{ } \beta
$$
since $\beta \cdot \mathbf{1} \in \tfrac{1}{L!} P(\mathbf{a}_{L!})$.  So
$$
\textnormal{Arn}(\mathbf{a}_{\bullet})  \textnormal{ } \leq  \textnormal{ } \beta  \textnormal{ } <  \textnormal{ } \lambda + \epsilon.
$$

For (3), one has
\begin{eqnarray*}
e(\mathbf{a}_{\bullet}) &=& k! \lim_{L \rightarrow \infty} \textnormal{Vol} \big( \tfrac{1}{L!} P(\mathbf{a}_{L!})^c \big) \\
&\leq& k! \lim_{L \rightarrow \infty} \textnormal{Vol} \big( (\cup_{n=1}^L \tfrac{1}{n} P(\mathbf{a}_{n}))^c \big) \\
&=& k! \cdot \textnormal{Vol} \big( P ( \mathbf{a}_{\bullet})^c \big). 
\end{eqnarray*}
To pass from the second to the third lines above, one uses that the sequence of sets
$$
\mathbf{R}_{\geq 0}^k \setminus \cup_{n=1}^L \tfrac{1}{n} P(\mathbf{a}_{n})
$$
is a nested, decreasing sequence of sets of finite volume, whose common intersection is $\mathbf{R}_{\geq 0}^k \setminus P(\mathbf{a}_{\bullet})$.  $\Box$

\begin{corollary}
The values of $\textnormal{ord}_0 (\mathbf{a}_{\bullet})$, $\textnormal{Arn} (\mathbf{a}_{\bullet})$, and $e(\mathbf{a}_{\bullet})$ only depend on the closure $\overline{P (\mathbf{a}_{\bullet})}$.
\end{corollary}
$\textbf{Proof}$.  It is clear from the previous proposition that the first two invariants depend only on $\overline{P (\mathbf{a}_{\bullet})}$.  Since $P(\mathbf{a}_{\bullet})$ and $\overline{ P (\mathbf{a}_{\bullet}) }$ differ only by a set contained in $\partial  P (\mathbf{a}_{\bullet})$ which has measure $0$ (the boundary of a convex set in $\mathbf{R}^n$ has measure $0$: see section 2.6 in [8]), it follows that
$$
e(\mathbf{a}_{\bullet}) = \textnormal{Vol} \big( P(\mathbf{a}_{\bullet}) \big) = \textnormal{Vol} \big( \overline{P(\mathbf{a}_{\bullet})} \big). \textnormal{ } \Box 
$$

In the next section, we will compute asymptotic invariants by working directly with closed, convex sets absorbing $\mathbf{R}_{\geq 0}^k$ under addition, without explicitly keeping track of the underlying graded systems of ideals.  Thanks to the previous corollary, the following construction suffices.

\begin{proposition}
Given a nonempty closed convex set $P \subseteq \mathbf{R}_{\geq 0}^k $ absorbing $\mathbf{R}_{\geq 0}^k$ under addition, there exists an $\mathbf{N}$-graded system of monomial ideals $\mathbf{a}_{\bullet}$ so that $\overline{ P(\mathbf{a}_{\bullet})} = P$.
\end{proposition}
\textbf{Proof}  Let $ \mathbf{a}_m$ be the monomial ideal defined by 
$$
\mathbf{a}_m := \langle \{ x^v | v \in mP \cap \mathbf{Z}^k \} \rangle
$$
We claim these form an $\mathbf{N}$-graded system, and moreover that $\overline{ P(\mathbf{a}_{\bullet})} = P$.  

To check that $\mathbf{a}_m \mathbf{a}_n \subseteq \mathbf{a}_{m+n}$, it suffices to show that if $v$ and $w$ are exponent vectors for monomials in $\mathbf{a}_n$ and $\mathbf{a}_m$ respectively, then $v+w$ is an exponent vector for a monomial in $\mathbf{a}_{m+n}$.  This reduces to showing that $mP + nP \subseteq (m+n)P$, which follows from convexity: if $v = np_1$ and $w = mp_2$ for $p_i \in P$, 
$$
v + w = mp_1 + np_2 = \big( m+n \big) \big( \tfrac{m}{m+n} p_1 + \tfrac{n}{m+n} p_2 \big) \in \big( m+n \big) P.
$$

We now establish that  $\overline{ P(\mathbf{a}_{\bullet})} = P$.  Clearly $P(\mathbf{a}_{\bullet}) \subseteq P$.  Since the condition that $P$ absorbs $\mathbf{R}_{\geq 0}^k$ under addition forces its interior $P^0$ to be dense in $P$, it suffices to show that $P^0 \subseteq P(\mathbf{a}_{\bullet})$.  Let $v \in P^0$ and select $w \in P^0 \cap \mathbf{Q}^k$ such that $v-w \in \mathbf{R}_{\geq 0}^k$.  Select $j \in \mathbf{N}$ clearing the denominators of the entries of $w$.  Then $jw \in \mathbf{Z}^k \cap jP$, so $w \in \frac{1}{j} P(\mathbf{a}_j)$.  Since $ \frac{1}{j} P(\mathbf{a}_j)$ absorbs $\mathbf{R}_{\geq 0}^k$ under addition, $v \in  \frac{1}{j} P(\mathbf{a}_j) \subseteq P(\mathbf{a}_{\bullet})$.  $\Box$

\section[Examples]{Behavior of Cones and Asymptotic Invariants}

In this section, we show that any closed convex cone can occur as the nef cone of a graded system of ideals.  From the construction, it follows that the asymptotic invariants associated to a graded system of ideals need not be even piecewise $C^{\infty}$ away from the nef cone, where they are identically zero.  We give another construction showing that even graded systems of ideals with almost everywhere locally analytic nef cones need not have asymptotic invariants that are $C^{\infty}$ on a dense open set in their domains.

We will rely on a few facts about convex cones which we collect here.  Let $C$ be a closed convex cone lying in a finite-dimensional real vector space $V$.
\begin{enumerate}
  \item If $C \neq V$ then there exists a hyperplane $L \subset V$ so that $C$ lies to one side of $L$.  (Theorem 1, Chapter 1 of [4])  \label{hyper}
  \item If \textnormal{Int}(C) $\neq \emptyset$, then $C = \overline{ \textnormal{Int}(C)}$.  (Proposition 14, Chapter 2 of [4]) \label{Int}  
    \item  If \textnormal{Int}(C) $\neq \emptyset$, by using the previous fact it is simple to show that C is the closed convex cone spanned by its integer points: it suffices to show that the integer points generate the Int($C$).  But the integer points generate the rational points of Int$(C)$, and since these are dense their convex hull contains all of Int$(C)$. \label{Z points}
      \item Suppose $C \subset \mathbf{R}^{n+1}$ with coordinates $x_1, \ldots, x_n, y$.  Let $L$ be the hyperplane with coordinates $x_1, \ldots, x_n$ and let $\pi: V \rightarrow L$ denote projection.  Suppose that $C$ is given by $y \geq f(x_1, \ldots, x_n)$ where $(x_1, \ldots, x_n)$ ranges over the convex set $\pi(C) \subseteq L$.  Then $f$ is convex on $\pi(C)$ and continuous on the interior of $\pi(C)$.  (See Propositions 1 and 23 in Chapter 3 of [4]) \label{conefunction}
            \item There exist convex functions on $\mathbf{R}^{2}$ that fail to be differentiable on any open set $U \subseteq \mathbf{R}^2$.
\end{enumerate}

Although the existence of convex functions on $\mathbf{R}^n$ that fail to be differentiable on any open set seems to be ``well-known,'' we have been unable to find a clear proof in the literature on convex functions.  In the interests of completeness, we give a simple example of such when $n=2$ in the Appendix.

\medskip

The following theorem shows that there are no restrictions besides convexity on the nef cone of a multigraded system of ideals.

\begin{theorem}
Given a closed convex cone $C$ contained in a finite-dimensional real vector space $V$, with nonempty interior, and with the vertex of $C$ lying at the origin of $V$, there is a graded system of ideals having $C$ as its nef cone.
\end{theorem}
$\textbf{Proof}$  We first concentrate on the case where $C \subset \mathbf{R}^{n+1}$ is given by $y \geq f(x_1, \ldots, x_n)$, where $f$ is defined and nonnegative at all $(x_1, \ldots, x_n) \in \mathbf{R}^n$.  After proving the theorem in this case, we will show how to find coordinates on $V$ putting $C$ in this form.

According to (4), $f$ is convex and continuous on $\mathbf{R}^n$.  Since $f$ defines a closed convex cone, moreover, $f$ is positively homogeneous:
$$
f( \lambda \cdot (x_1, \ldots, x_n) ) = \lambda \cdot f(x_1, \ldots, x_n), \textnormal{ for } \lambda \in \mathbf{R}_{\geq 0}.
$$
We set  $\lceil x \rceil$ to be the least integer $n \geq x$.  The properties of $f$ and of $\lceil \cdot \rceil$ imply that the map $\phi: \mathbf{Z}^{n+1} \rightarrow \mathbf{Z}$ defined by 
$$
\phi(x_1, \ldots, x_n, y) := \lceil f(x_1, \ldots, x_n) - y \rceil
$$
is subadditive and positively homogeneous on $\mathbf{Z}^{n+1}$.

Quite generally, if $\phi$ is any such mapping and $I_\bullet$ is a $\mathbf{Z}$-graded system of ideals with $I_{m+1} \subseteq I_m$ for all $m \in \mathbf{Z}$, then we may define a $\mathbf{Z}^{n+1}$-graded system $\mathbf{a}_\bullet$ by
$$
\mathbf{a}_v := I_{\phi(v)}.
$$
We work in the ring $R = \mathbf{C}[ z_1, z_2 ]$ and run the construction outlined above with the graded system $I_m = (z_1, z_2)^m$.  Note that we are taking $I_m =R$ if $m \leq 0$.

We thus obtain a $\mathbf{Z}^{n+1}$-graded system of ideals $\mathbf{a}_\bullet$ so that if $\overline{v} = (x_1, \ldots, x_n, y)$,
$$
\mathbf{a}_{\overline{v} } = (z_1, z_2)^{ \lceil f(x_1, \ldots, x_n) -y \rceil }.
$$
Now $\overline{v}$ satisfies $\mathbf{a}_{ \overline{v}} = R$ precisely if $\lceil f(x_1, \ldots, x_n) - y \rceil \leq 0$.  But this conditions simply says that the integer point $\overline{v}$ lies in $C$.  Thus $\textnormal{Nef}(\mathbf{a}_\bullet) \subseteq C$.  However, (3) implies that in fact $\textnormal{Nef}(\mathbf{a}_\bullet) = C$.

We will denote this graded system of ideals by $\mathbf{a}_\bullet = \mathbf{m}_{\bullet}^C$ to indicate its dependence on $C$.

We now show how to obtain the system of coordinates on $V$ in which $C$ has the desired form.  We assume $C \subsetneq V$, otherwise the statement of the theorem is trivial.  Therefore by (1) we may fix a hyperplane $L$ through $0$ with $C$ to one side of $V$.  Let $e_1, \ldots e_n$ be a basis for $L$ and let $v$ be in the interior of $C$ so that $e_1, \ldots, e_n, v$ form a basis for $V$.  Using the identification $V \rightarrow \mathbf{R}^{n+1}$ afforded by this basis, one obtains a projection $\pi: V \rightarrow L$ and an isomorphism $\mathbf{Z}^{n+1} \rightarrow \sum \mathbf{Z}e_i + \mathbf{Z}v$.

We further observe that since $v \in \textnormal{Int}(C)$, $\pi(C) = L$, so above any $\overline{x} \in L$, there is a point $(\overline{x},y) \in C$ for some $y \geq 0$.  We define $f: L \rightarrow \mathbf{R}$ by
$$
f(x_1, \ldots, x_n) := \textnormal{inf} \{ y | \sum_{i=1}^n x_i e_i + yv \in C \}. 
$$
Then $f$ is nonnegative convex, continuous, and positively homogeneous on $\mathbf{R}^n$.  $\Box$

\medskip
This gives an example of a multigraded system which has the entire algebraic Neron-Severi space as its effective cone, as mentioned in Section 2.

The construction in the proof of the theorem also shows that the asymptotic invariants of a general multigraded system of ideals are exotic functions on Big$(\mathbf{a}_{\bullet})$.

\begin{corollary}
There exist graded systems $\mathbf{a}_{\bullet}$ whose asymptotic invariants are non-$C^{\infty}$ functions on a dense open subset of $\textnormal{Big}(\mathbf{a}_{\bullet}) \setminus \textnormal{Nef} (\mathbf{a}_{\bullet})$.
\end{corollary}
\textbf{Proof}.  Keeping notation as in the previous theorem, let $\mathbf{a}_{\bullet} = \mathbf{m}_{\bullet}^C$, for $C$ a closed convex cone whose defining equation $f$ is nowhere $C^{\infty}$ in the variables $x_1, \ldots, x_n$.  As noted in (5), these exist at least for $n=2$.  Let $(x_1, \ldots, x_n, y) \in \textnormal{Big}(\mathbf{a}_{\bullet}) \setminus \textnormal{Nef} (\mathbf{a}_{\bullet})$ be an integer-valued point.
\begin{eqnarray*}
\textnormal{ord}_0(\mathbf{m}_{\bullet}^C, (x_1, \ldots, x_n,y)) &=& \lim_k \frac{\textnormal{ord}_0 (\mathbf{m}^{ \lceil f(kx_1, \ldots, kx_n) - ky \rceil } ) }{k} \\
  &=& \lim_k \Big (  \frac{\lceil k(f(x_1, \ldots, x_n) - y ) \rceil }{k} \Big ) \textnormal{ord}_0 (\mathbf{m})  \\
  &=& f(x_1, \ldots, x_n) - y.
\end{eqnarray*}
\begin{eqnarray*}
\textnormal{Arn}(\mathbf{m}_{\bullet}^C, (x_1, \ldots, x_n, y)) &=& \lim_k \frac { \textnormal{Arn} (\mathbf{m}^{ \lceil f(kx_1, \ldots kx_n) - ky \rceil})}{k} \\
&=&  \lim_k  (\frac{ \lceil k ( f(x_1, \ldots, x_n) - y) \rceil }{k}) \textnormal{Arn} (\mathbf{m}) \\
&=& \frac{f(x_1, \ldots, x_n) - y}{L}.
\end{eqnarray*}
\begin{eqnarray*}
e_{\mathbf{m}}( \mathbf{m}_{\bullet}^C, (x_1, \ldots, x_n,y)) &=&
\lim_k \frac{ e_{\mathbf{m}} ( \mathbf{m}^{ \lceil f(kx_1, \ldots, kx_n) - ky \rceil })}{k^L} \\
&=& \lim_k (\frac{ \lceil f(kx_1, \ldots kx_n) - ky \rceil }{k} )^L e_{\mathbf{\mathbf{m}}}(\mathbf{m}) \\
&=& ( f(x_1, \ldots x_n) - y)^L.
\end{eqnarray*}

By homogeneity these formulas hold for rational-valued points in $ \overline{\textnormal{Eff}}(\mathbf{a}_{\bullet}) \setminus \textnormal{Nef} (\mathbf{a}_{\bullet})$, and hence by continuity for all real-valued points as well.  The lemma now follows from (5).  $\Box$

\medskip
We observe that in the previous example the properties of the asymptotic invariants are essentially described by the geometry of the cone $C$.  Since it is known that the nef cone of a projective variety is almost everywhere given by analytic inequalities (see [1]), we wish to give an example of a multigraded system with locally analytic nef cone but non-differentiable order function.  

\begin{theorem}
There exists a graded system $\{ \mathbf{a}_v \}_{v \in \mathbf{Z}^2}$ so that the boundary of $\textnormal{Nef}(\mathbf{a}_{\bullet})$ is given by linear functions, but for which $\textnormal{ord}_0$ fails to be differentiable on an open set in $\textnormal{NS}_{\mathbf{R}}(\mathbf{a}_{\bullet})$.
\end{theorem}

In the example we will construct, Nef$(\mathbf{a}_{\bullet})$ is the third quadrant and $\overline{\textnormal{Eff}} (\mathbf{a}_{\bullet}) = \mathbf{R}^2$.  We'll observe in Remark 4 that one modify this construction to obtain a graded system whose effective cone lies to one side of a line.

$\textbf{Proof}$.  Fix two real, nonnegative, decreasing convex functions $f$ and $g$ on $\mathbf{R}_{\geq 0}$ and let $P$ and $Q$ denote, respectively, the sets in $(\mathbf{R}_{\geq 0})^2$ above the graphs of $f$ and $g$.  Let $\mathbf{a}_{\bullet}$ and $\mathbf{b}_{\bullet}$ be the $\mathbf{N}$-graded systems of monomial ideals on $\mathbf{C}^2$ with $P(\mathbf{a}_{\bullet}) = P$ and $P ( \mathbf{b}_{\bullet} ) = Q$, described in the proof of Proposition 2.  We extend these to $\mathbf{Z}$-graded systems by setting
$$
a_n = b_n = \mathbf{C} [ x_1, x_2 ], \textnormal{ for } n \textnormal{ } \leq \textnormal{ }0 
$$
We obtain a $\mathbf{Z}^2$-graded system of ideals $ \mathbf{c}_{\bullet} := pr_1^{*}(\mathbf{a}) \cap pr_2^{*} (\mathbf{b})$ defined by
$$
\mathbf{c}_{(m,n)} = \mathbf{a}_m \cap \mathbf{b}_n.
$$
$\mbox{Nef}(\mathbf{c}_{\bullet})$ is the third quadrant, so it is described by linear functions.  However, we will show that some of the asymptotic invariants for this graded system are not differentiable on open sets in the first quadrant, for appropriate choices of the functions $f$ and $g$.  Let $g(x) = 1 - \frac{x}{2}$ and let $f$ be a convex function that fails to be differentiable on a dense set of points: we will add a dense set of ``kinks'' by choosing
$$
f(x) = (-2x + 2) + \sum_i \textnormal{ max } \{ 0, \frac{\epsilon_i - x}{n_i} \}$$ 
where $\epsilon_i$ are dense in $[0,1]$ and the $n_i$ are chosen so that $\sum_i \tfrac{\epsilon_i}{n_i} < 1$.

We concentrate on $\textnormal{ord}({\mathbf{c}}_{\bullet},(r,s))$ where $(r,s)$ is in the first quadrant.  We observe that in the first quadrant of $(m,n)$-space the ideal $\mathbf{c}_{(m,n)}$ is generated by all mononmials $x^v$ where $v$ is a lattice point in $mP \cap nQ$.  Therefore one has that
$$
\textnormal{ord}_0(\mathbf{c}_{\bullet}, (m,n)) = \mbox{inf} \{ x + y | (x,y) \in mP \cap nQ \}
$$
By homogeneity, these results  hold true for $m,n$ positive rational numbers, and hence by continuity for $r,s$ positive real numbers.

For $(r,s)$ in a sufficiently small neighborhood $U$ of $(1,1)$, $f_r$ and $g_s$ meet in a single point $p$. We claim that $p$ is the point at which the infimum giving $\textnormal{ord}_0$ is attained: since $f_r$ decreases faster than a linear function with slope $-1$ and the slope of $g_s$ is $\tfrac{1}{2}$, the line of the form $x+y = C$ through $p$ passes under all other points in $rP \cap sQ$.  We write $p = (x_{r,s}, s - \frac{ x_{r,s} }{2})$.

So on the open set $U$, we have that $\textnormal{ord}_0 (r,s) = s + \tfrac{x_{r,s}}{2}$.  To establish the theorem, it is therefore enough to show $x_{r,s}$ is not differentiable on $U$. We view $r$ as fixed. 

Since $f_r$ fails to be differentiable on a dense set, on a dense set of $s-$values the graph of $g_s$ meets that of $f_r$ at over an $x-$value where $f_r$ fails to be differentiable.  We claim that $x_{r,s}$ is not differentiable with respect to $s$ at those points $(r, s_0)$.

Since $f_r$ has a kink at $(r, s_0)$, there exists a line $\ell$ of slope $\sigma < -2 $ through $(r, s_0)$ so that for $s > s_0$
$$
\frac{ f_r (s) - f_r (s_0 }{s-s_0} > \sigma
$$
and for $s < s_0$
$$
\frac{ f_r (s) - f_r (s_0 }{s-s_0} < \sigma - \epsilon
$$
If $y_s$ denotes the intersection of $\ell$ and $g_s$, then one finds
$$
\frac{y_s - y_{s'} }{s - s'} = \frac{1}{\sigma + \tfrac{1}{2} }.
$$
Comparing the intersection points of $g_s$ and $f_r$ to those of $\ell$ and $g_s$, one sees that for $s > s_0$
$$
\frac{ x_{r,s} - x_{r, s_0} }{s-s_0} < \frac{1}{ \sigma + \tfrac{1}{2} }
$$
but for $s < s_0$
$$
\frac{ x_{r,s} - x_{r, s_0} }{s-s_0} > \frac{1}{ \sigma + \tfrac{1}{2} - \epsilon }
$$
So the left and right difference quotients for $x_{r,s}$ are bounded away from one another. $\Box$

\begin{remark}
\textnormal{By truncating the graded system by the semigroup of points } $\{ (r,s) : s \geq \epsilon r \}$ \textnormal{ for small, positive } $\epsilon$, \textnormal{one preserves the region where the variation of} $\textnormal{ord}_0$ \textnormal{ is non-} $C^{\infty},$ \textnormal{while arranging that} $\overline{\textnormal{Eff}}(\mathbf{a}_{\bullet})$ \textnormal{ lies to one side of a line.}
\end{remark}

\section{Appendix}

\begin{proposition}
There exist convex functions on $\mathbf{R}^2$ that are nowhere differentiable on a dense open set $\mathbf{R}^2$.
\end{proposition}

\textbf{Proof} Let $R'$ denote the region in $\mathbf{R}_{\geq 0}^2$ under the graph of the concave function on $[0,1]$ defined by
$$
f(x) = \sum_{i} \textnormal{ min } \{ \epsilon_i, \frac{\epsilon_i}{1 - x_i} (x-1) \}
$$
where $\sum_i \epsilon_i$ converges and the $x_i$ form a dense set on $[0,1]$. Reflecting $R'$ about the $x-$ and $y-$axes, we obtain a convex set with a nowhere differentiable boundary.

Let $C \subset \mathbf{R}^3$ denote the cone in $\mathbf{R}^3$ over $R$ viewed as a subset of $\mathbf{R^2} \times 1 \subset \mathbf{R}^3$, and let $F$ be the function on $\mathbf{R}^2 \times 0$ whose epigraph is $C$.

We claim $F$ is nowhere differentiable on $\mathbf{R}^2$. By the construction of $R$ it suffices to show this on the first quadrant. Let $(s_0, t_0, 0)$ lie under a point $p$ on the boundary of $R$ where $R$ fails to be differentiable. Then there exist lines $l_1, l_2$ with distinct slopes $\lambda_1 < \lambda_2$ and $y-$intercepts $d_1 > d_2$, passing through $p$, lying under the graph of $f$ if $s < s_0$ and lying over the graph of $f$ if $s > s_0$. We compare $F$ with the functions 
$$
F_i (s,t) = \frac{t_0 - \lambda_i s_0}{d_i}
$$
whose epigraphs are cones over those lines. For $s < s_0,$
\begin{eqnarray*}
\frac{F(s,t_0) - F(s_0, t_0)}{s-s_0} &\geq& \big( \frac{t_0 - \lambda_1 s_0}{d_1} - \frac{t_0 - \lambda_1 s_0}{d_1} \big) /(s-s_0) \\
&=& - \frac{ \lambda_1 }{d_1} \\
&=& - \frac{t_0}{d_1 s_0} + \frac{1}{s_0}.
\end{eqnarray*}
Similarly, for $s > s_0$, 
$$
\frac{F(s,t_0) - F(s_0, t_0)}{s-s_0} \leq - \frac{t_0}{d_2 s_0} + \frac{1}{s_0}
$$
and one checks that the right and left hand difference quotients are bounded away from one another. 

Since $F$ is positively homogeneous, this shows that $F$ fails to be differentiable along any ray through points $(s_0, t_0)$. Since such rays are dense, there can be no open set in $\mathbf{R}^2$ on which $F$ is differentiable.  $\Box$

\end{document}